\documentclass[10pt,twocolumn,twoside]{IEEEtran}

\usepackage{mathbbold}
\usepackage{color}
\usepackage{amsfonts}
\usepackage{bbm}
\usepackage{mathrsfs}
\usepackage{amsmath}
\usepackage{graphics}
\usepackage{epsfig}
\usepackage{amssymb}

\newtheorem{rema}{Remark}
\newtheorem{theo}{Theorem}

\newtheorem{lemm}{Lemma}
\newtheorem{coro}{Corollary}
\newtheorem{exam}{Example}

\newtheorem{prop}{Proposition}
\begin{document}
\title{\LARGE \bf Random Asynchronous Iterations in Distributed\\ Coordination Algorithms}
\author{Yao Chen, Weiguo Xia, Ming Cao, and Jinhu L\"{u}
\thanks{Y. Chen is with the Department of Computer Science, Southwestern
        University of Finance and Economics, Chengdu 611130, China (e-mail: chenyao07@gmail.com). W. Xia is with the School of Control Science and Engineering, Dalian University of Technology, Dalian 116024, China  (e-mail: wgxiaseu@dlut.edu.cn). M. Cao is with the Faculty of Science and Engineering, ENTEG,
        University of Groningen, Groningen 9747 AG, the Netherlands (e-mail: m.cao@rug.nl). Jinhu L\"{u} is with the School of Automation Science and Electrical Engineering, State Key Laboratory of Software Development Environment, Beijing Advanced Innovation Center for Big Data and Brain Computing, Beihang University, Beijing 100191, China, and the Academy of         Mathematics and Systems Science, Chinese Academy of Sciences, Beijing 100190, China  (e-mail: jhlu@iss.ac.cn).}
}


\maketitle

\begin{abstract}
Distributed coordination algorithms (DCA) carry out information processing processes among a group of networked agents without centralized information fusion. Though it is well known that DCA characterized by an SIA (stochastic, indecomposable, aperiodic) matrix generate consensus asymptotically via synchronous iterations, the dynamics of DCA with asynchronous iterations have not been studied extensively, especially when viewed as stochastic processes. This paper aims to show that for any given irreducible stochastic matrix, even non-SIA, the corresponding DCA lead to consensus successfully via random asynchronous iterations under a wide range of conditions on the transition probability. Particularly, the transition probability is neither required to be independent and identically distributed, nor characterized by a Markov chain.
\end{abstract}

\section{Introduction}

Distributed coordination algorithms (DCA) use local information of a group of networked agents to generate specific collective behaviors. Nowadays, DCA have been used not only to explain social or economic phenomena \cite{JASA-DeGroot}, but also to solve practical engineering problems \cite{I3EAC-Nedic-Olshevsky-2014}, \cite{EURASIP-XGang}. A typical DCA generates aligned collective motion, usually referred to as consensus (or heading synchronization, see \cite{PRL-Vicsek}, \cite{Yu-2nd}).

When agents update their states using DCA, they may do so synchronously or asynchronously. More precisely, by synchronous updating we mean all the agents update their states at the same time repeatedly. This requires the availability of a global clock or a set of identical local clocks, which is a stringent requirement in practice (see, \cite{JCSS-Dolev}). By asynchronous updating we mean each agent has an independent local clock, according to which it updates its own state without paying attention to when the other agents update. Only in this way DCA can be constructed in a fully distributed manner \cite{MingCao_AgreeAsynchrous}. It has been recently reported that the synchronous and asynchronous implementation of DCA may lead to dramatically different asymptotic collective behaviors; for example, some DCA may converge under synchronous updating but diverge under asynchronous updating \cite{I3EAC-Xia-Cao}.

The possible significant differences of DCA in deterministic and stochastic settings are reflected in the fundamental differences in their corresponding analytical tools: the determination of convergence in DCA, even for a pair of stochastic matrices, has been proved to be NP-hard in a deterministic setting \cite{SIAMCO-Blondel}; in sharp comparison, the convergence for DCA with finite Markovian random switching modes can be determined by using classic LMI-based techniques \cite{MCSS-Boukas}, which have been proved to obtain desired solutions with high efficiency. Therefore, the analysis in stochastic settings often generates less conservative sufficient conditions for consensus: In 2007, Porfiri and Stilwell gave some sufficient conditions for consensus over random weighted directed graphs generated by independent and identically distributed (i.i.d) random variables \cite{I3EAC-Porfiri}. In 2013, You {\em et al.} established some necessary and sufficient conditions for consensus based on the assumption of Markovian switching topologies \cite{Automatica-You}; Matei {\em et al.} gave some sufficient conditions for the linear consensus problem under Markovian random graphs. Note that the network topologies in these works are either Markovian or generated by i.i.d random variables. When asynchronous updating in DCA is not generated by a Markovian chain or independent random variables, the analysis becomes much more challenging due to the limitations of the existing analytical methods.


In the deterministic setting, the topological condition for consensus of DCA with asynchronous iteration is generally very restrictive. In 2014, Xia and Cao proved that for any given scrambling matrix, the corresponding DCA reach consensus for any asynchronous iteration \cite{I3EAC-Xia-Cao}. It should be noted that any pair of nodes in the graph of a scrambling matrix share a common neighbor, which makes such a graph densely connected. In this paper, we will investigate the asynchronous iterations of DCA in the stochastic setting with the aim of relaxing the restrictive constraint on topological structures. To realize this purpose, we will transform the consensus problem of random asynchronous DCA to the random walk problem along a labelled directed cycle, and propose a graphical method to analyze the convergence of random asynchronous DCA. Specifically, the contributions of this paper can be summarized as follows:
\begin{itemize}
\item[a).] The convergence of asynchronous iterations of DCA is investigated in the stochastic setting for the first time. The obtained results only require the graph of the given matrix of DCA to be connected. Compared with the related results in the deterministic setting, the matrix is not required to be SIA or scrambling;
\item[b).] The critical conditions for consensus in the traditional stochastic setting, such as i.i.d or Markovian switching, and positive diagonal entries of switching matrices are not needed in our main result any more.
\end{itemize}

The rest of the paper is organized as follows: Section \ref{se:2} gives some preliminaries on graph theory and formulates the problem of random asynchronous updating for DCA; Section \ref{se:3} presents the main results and related discussions; Section \ref{se:4} provides the skeleton of the technical proof; Section \ref{se:5} gives some numerical examples and Section \ref{se:6} concludes this paper. More details of the proof of the main theorem, examples and corollaries can be found in the appendix.

\section{Preliminaries and Problem Formulation}\label{se:2}

In this section, we will give some preliminaries on graph theory and asynchronous iterations of DCA in two subsections, respectively.

\subsection{Preliminaries}
A graph $\mathcal{G}=(V,\mathcal{E})$ is composed of two sets, where $V$ is the set of nodes and $\mathcal{E}\subseteq V\times V$ is the set of edges. A path of $\mathcal{G}$ is composed of a sequence of distinct nodes $i_1, i_2, \cdots i_k$ which satisfy $(i_j, i_{j+1})\in\mathcal{E}$ for any $1\leq j\leq k-1$. $\mathcal{G}$ is called rooted if there exists a node $r\in V$ such that for any $j\neq r$, there is a path from $r$ to $j$, where $r$ is called a root of $\mathcal{G}$.
$\mathcal{G}$ is called strongly connected if there exists a path from any node $i\in V$ to any node $j\in V$ ($i\neq j$). The collection of all the roots of graph $\mathcal{G}$ is defined as $\mathbbm{r}(\mathcal{G})$. If there exists an edge from node $i$ to itself, we call such an edge a self-loop. Specifically, if $\mathbbm{r}(\mathcal{G})=\{r\}$, we say node $r$ is the unique root of graph $\mathcal{G}$.

In order to regroup the nodes in a graph, we assign a positive integer to each node of a graph, and these integers are called labels of the nodes. A graph $\mathcal{G}=(V_l,\mathcal{E})$ is said to be a labelled graph if each node $i\in V_l$ is assigned with an integer label ${\rm label}(i)$. It should be noted that there may exist two nodes $i,j\in V_l$ with identical labels, i.e., ${\rm label}(i)={\rm label}(j)$. The following proposition tells us the fact: any strongly connected graph can be mapped to a labelled directed cycle.

\begin{prop}\label{prop:graphcycle}
Given any strongly connected graph $\mathcal{G}=(V,\mathcal{E})$ with $N$ nodes, it can be mapped to a labelled directed cycle $\mathscr{C}$ with length $l\leq N(N-1)$ that contains all the nodes of $V$.
\end{prop}

{\em Proof: } Due to the strong connectivity of graph $\mathcal{G}$, there exists a path from $k$ to $k+1$ ($1\leq k\leq N-1$) with length no more than $N-1$, and there also exists a path from node $N$ to $1$ with length no more than $N-1$. Concatenate these $N$ paths and one generates a cycle which contains all the nodes of $V$. \hfill $\square$

\begin{figure}[htbp]
\begin{center}
\resizebox{7.2cm}{2.8cm}{\includegraphics{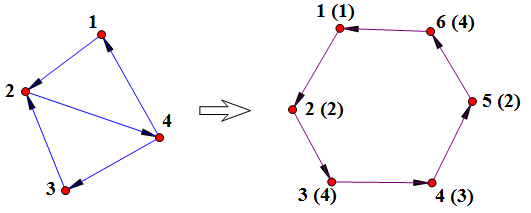}}
\end{center}
\caption{Transform a strongly connected graph to a labelled directed cycle. }\label{fig:cycle}
\end{figure}

As shown in Fig. \ref{fig:cycle}, a strongly connected graph with $4$ nodes can be transformed to a labelled directed cycle with $6$ nodes, in which $V_l=\{1,2,3,4,5,6\}$ and the label of each node $i\in V_l$ is given in the bracket. In particular, nodes $2$ and $5$ share the same label $2$, nodes $6$ and $3$ share the same label $4$.

A matrix $A=(a_{ij})_{i,j=1}^N$ is called stochastic (or row stochastic) if $a_{ij}\geq 0$ ($\forall\, i,j\in V$) and $\sum_{j=1}^N a_{ij}=1$ ($\forall i\in V$). The graph $\mathcal{G}(A)=(V,\mathcal{E})$ corresponding to $A$ is defined by: $V=\{1,2,\cdots, N\}$, and $(j,i)\in\mathcal{E}$ if and only if $a_{ij}>0$. Based on $\mathcal{G}(A)$ we define $\mathbbm{r}(A)=\mathbbm{r}(\mathcal{G}(A))$ and
\begin{eqnarray}
\mathcal{N}(A,\nu) = \{k: a_{\nu k}>0,\,\, k\in V\} \nonumber
\end{eqnarray}
as the neighbor of node $\nu$ in $\mathcal{G}(A)$. The ergodic coefficient of a stochastic matrix $A$ is defined by
\begin{eqnarray}
\lambda(A) = 1-\min_{i\neq j}\sum_{k=1}^N\min(a_{ik},a_{jk}). \nonumber
\end{eqnarray}
A stochastic matrix $A$ is called scrambling if $\lambda(A)<1$, and is called SIA if $\lim_{k\rightarrow\infty}A^k=\mathbbold{1}\xi^T$ for some $\xi\in\mathbb{R}^N$, where $\mathbbold{1}\in\mathbb{R}^N$ is a vector with each entry being $1$.

Given two stochastic matrices $A=(a_{ij})_{i,j=1}^N,B=(b_{ij})_{i,j=1}^N\in\mathbb{R}^{N\times N}$, we say $A$ and $B$ are of the same type if ${\rm sgn}(a_{ij})={\rm sgn}(b_{ij})$ holds for any $i,j=1,2,\cdots N$, denoted by $A\sim B$, where ${\rm sgn}(\cdot)$ is the sign function.

Given a vector $x=(x_1,x_2,\cdots x_N)^T\in \mathbb{R}^N$, we define
\begin{eqnarray}
\Delta(x) = \max_{i=1}^N x_i - \min_{i=1}^{N}x_i \nonumber
\end{eqnarray}
as the maximal discrepancy of vector $x$. Given a stochastic matrix $A$ and $y=Ax$ with $y\in\mathbb{R}^N$, it holds that \cite{CaoM-SIAM2}
\begin{eqnarray}
\Delta(y)\leq \lambda(A)\Delta(x). \nonumber
\end{eqnarray}
Specifically, the ergodic coefficient $\lambda(\cdot)$ has the following important properties.
\begin{itemize}
\item[a)] Given any stochastic matrices $A_1, A_2\in\mathbb{R}^{N\times N}$, it holds that \cite{MingCao_AgreeAsynchrous}
    \begin{eqnarray}
    \lambda(A_1A_2)\leq \lambda(A_1)\lambda(A_2). \nonumber
    \end{eqnarray}
\item[b)] Given any stochastic matrix $A\in\mathbb{R}^{N\times N}$, it holds that \cite{Springer-Senate}
    \begin{eqnarray}
    \lambda(A) = \sup_{\Delta(x)=1}\Delta(Ax). \nonumber
    \end{eqnarray}
\end{itemize}

A matrix $A=(a_{ij})_{i,j=1}^N$ is called column stochastic if $A^T$ is row stochastic. In this paper, we use a column stochastic matrix to denote a Markovian chain. The entry $a_{ij}$ in a column stochastic matrix represents the transition probability from node $j$ to node $i$. Without specific declaration, for any stochastic matrix in this paper, we mean it is row stochastic.

\subsection{Problem Formulation}

DCA characterizes the evolution of states in a network of agents via local interaction. A typical one of DCA is the following linear averaging protocol
\begin{eqnarray}
x_i(k+1)= \left\{\begin{array}{ll}
\sum_{j=1}^N a_{ij}x_j(k), &\mbox{if agent } i \mbox{ updates } \nonumber\\
&  \mbox{at time intstant } k; \nonumber\\
x_i(k), &\mbox{if agent } i \mbox{ does not} \nonumber\\
&  \mbox{update at time instant } k, \nonumber
\end{array}\right.
\end{eqnarray}
where $x_i(k)$ is the state of agent $i$ at time $k$, $a_{ij}$ denotes the coupling coefficient between agent $i$ and $j$, $a_{ij}\geq 0$ for any $i,j\in V$, and $\sum_{j=1}^N a_{ij}=1$ for any $i\in V$. Define $A=(a_{ij})_{i,j=1}^N\in \mathbb{R}^{N\times N}$ as the corresponding stochastic coupling matrix.

If each agent $i\in V$ updates at any time $k$, then the above DCA transforms to
\begin{eqnarray}
x(k+1)=  Ax(k),\quad k=1,2,\cdots,  \label{eq:synchronous}
\end{eqnarray}
where $$x(k)=(x_1(k),x_2(k),\cdots, x_N(k))^T$$ represents the state vector for all the individuals. We call DCA (\ref{eq:synchronous}) the synchronous DCA.

However, synchronous DCA is difficult to be implemented in practice since the local clocks associated with all the nodes are generally nonidentical (see \cite{JCSS-Dolev}, the impossibility of clock synchronization). In the case of asynchronous iteration, only part of the agents $\sigma_k\in 2^V$ ($2^V$ is the set of all subsets of $V$) update their states at each time instant $k$, then one generates the following asynchronous DCA:
\begin{eqnarray}
x(k+1)=  A_{\sigma_k}x(k), \label{eq:asynchronous}
\end{eqnarray}
where $A_{\sigma_k}$ is the asynchronous iteration matrix with the index set $\sigma_k\subseteq V$, whose definition is given by the following construction process: if $j\in\sigma_k$, then the $j$th row of $A_{\sigma_k}$ equals the $j$th row of $A$; if $j\notin\sigma_k$, then the $j$th row of $A_{\sigma_k}$ is the $j$th elementary vector $\mathbbm{e}_j$. Specifically, when only one agent updates at time $k$, $\sigma_k$ becomes a singleton and in this case we also use $\sigma_k$ to denote the agent $\sigma$ ($\sigma\in V$) that updates in the sequel. It follows that
$$A_\sigma\triangleq A_{\sigma_k}=(e_1, e_2, \cdots e_{\sigma-1}, a_\sigma^T, e_{\sigma+1}, \cdots, e_N)^T,$$
where $a_\sigma$ is the $\sigma$th row of $A$ and $e_j$ are elementary column vectors.

For example, given a stochastic matrix
\begin{eqnarray}
A = \left(\begin{array}{ccc}
0 & 1 & 0\\
0.2 & 0.8 & 0\\
0 & 0.7 & 0.3
\end{array}\right), \nonumber
\end{eqnarray}
$A_{2}$ and $A_{\{1,3\}} $ are defined as
\begin{eqnarray}
A_2 = \left(\begin{array}{ccc}
1 & 0 & 0\\
0.2 & 0.8 & 0\\
0 & 0 & 1
\end{array}\right), \quad
A_{\{1,3\}} = \left(\begin{array}{ccc}
0 & 1 & 0\\
0 & 1 & 0\\
0 & 0.7 & 0.3
\end{array}\right). \nonumber
\end{eqnarray}

Denote $\{A_{\sigma_k}\}_{k=1}^{\infty}$ as an asynchronous iteration sequence of matrix $A$. $\{A_{\sigma_k}\}_{k=1}^{\infty}$ is said to generate consensus in the deterministic setting if for any initial value $x(1)\in\mathbb{R}^N$, there exists $\xi\in\mathbb{R}$ such that
\begin{eqnarray}
\lim_{k\rightarrow\infty}x(k+1)= \mathbbold{1}\xi \nonumber
\end{eqnarray}
in DCA (\ref{eq:asynchronous}).

The convergence of an asynchronous iteration sequence of a stochastic matrix cannot  be determined by the SIA property of this matrix: the asynchronous implementation of an SIA matrix may not generate consensus, but the synchronous implementation of a non-SIA matrix may generate consensus (see, examples in \cite{I3EAC-Xia-Cao}).

In 2014, Xia and Cao proposed the following sufficient condition for consensus of asynchronous DCA (\ref{eq:asynchronous}) in the deterministic setting.

\begin{prop}\label{prop:xiacao}
If $A$ is scrambling and there exists $q>0$ such that $\bigcup_{k=j}^{j+q-1} \sigma_k= V$ for any $j\geq 1$, then the asynchronous iteration sequence $\{A_{\sigma_k}\}_{k=1}^{\infty}$ generates consensus.
\end{prop}

However, when we consider asynchronous iteration in the stochastic setting, i.e., $\{\sigma_k\}_{k\geq 1}$ are generated by random variables $\{\Xi_k\}_{k\geq 1}$, the dynamics of asynchronous DCA (\ref{eq:asynchronous}) will be quite different. A major challenge for the analysis is that the value of $\sigma_k$ and the corresponding transition probability may depend on its historic values: $\sigma_{k-1}$, $\sigma_{k-2}$, $\cdots$, $\sigma_1$. In this paper, we use $(\Omega, \mathscr{F}, \mathbb{P})$ to represent the probability space, where $\Omega$ is the sample space of $\sigma_k$, $\mathscr{F}$ is the $\sigma$-field, and $\mathbb{P}$ is the probability function.

When $\sigma_k\in 2^{V}$, the possible values of $A_{\sigma_k}$ is
\begin{eqnarray}
\mathfrak{A}= \{A_{\sigma}: \sigma\in 2^V\}, \nonumber
\end{eqnarray}
where $|\mathfrak{A}|=2^N$. When $\sigma_k\in V$, the possible values of $A_{\sigma_k}$ is
\begin{eqnarray}
\mathfrak{A}= \{A_{\sigma}: \sigma\in V\}, \nonumber
\end{eqnarray}
where $|\mathfrak{A}|=N$. Since the diagonal entries of $A$ may contain zeros, the set $\mathfrak{A}$ may also contain elements whose diagonal entries have some zeros. The existence of zero diagonal entries in the coupling matrices usually brings big challenges in analyzing the products of stochastic matrices (see, \cite{ChenY-TAC16}, \cite{I3EAC-Touri-Nedic}, and \cite{Automatica-Touri-Nedic}).

In the stochastic setting, the asynchronous DCA (\ref{eq:asynchronous}) is said to realize consensus almost surely if
\begin{eqnarray}
\lim_{k\rightarrow\infty}\mathbb{P}\left(\sum_{j=1}^N\left\|x_j(k)-\frac{1}{N}\sum_{i=1}^N x_i(k)\right\|^2\geq \varepsilon\right)=0 \label{eq:cond}
\end{eqnarray}
for any $\varepsilon>0$ and $x(1)\in\mathbb{R}^N$. If one defines the following projection matrix
\begin{eqnarray}
P = I - \frac{1}{N^2}\mathbbold{1}\mathbbold{1}^T, \nonumber
\end{eqnarray}
where $I$ is the identity matrix, then (\ref{eq:cond}) can be equivalently rewritten as
\begin{eqnarray}
\lim_{k\rightarrow\infty}\mathbb{P}\left(\|Px(k)\|
\geq\varepsilon\right)=0 \label{eq:equalform}
\end{eqnarray}
for any $\varepsilon>0$ and $x(1)\in\mathbb{R}^N$.

In the rest of this paper, for the simplicity of expression, we denote
\begin{eqnarray}
k_1:k_2 &=& \{k_1, k_1-1, \cdots, k_2\}, \nonumber\\
\sigma_{k_1:k_2} &=&  \{\sigma_{k_1}, \sigma_{k_1-1},\cdots, \sigma_{k_2}\}, \nonumber\\
A_{\sigma_{k_1:k_2}} &=& A_{\sigma_{k_1}}A_{\sigma_{k_1-1}}A_{\sigma_{k_1-2}}\cdots  A_{\sigma_{k_2+1}}A_{\sigma_{k_2}}\nonumber
\end{eqnarray}
for any $k_1\geq k_2\geq 1$. Particularly, it follows
\begin{eqnarray}
\sigma_{k:k} = \{\sigma_k\} \nonumber
\end{eqnarray}
and in the case of $\sigma_k\in V$, we do not distinguish the expressions of $\sigma_k$ and $\{\sigma_k\}$.

In the next section, we will give several sufficient conditions which guarantee almost sure consensus of asynchronous DCA (\ref{eq:asynchronous}).


\section{Conditions for Random Asynchronous Consensus}\label{se:3}

We briefly summarize the sufficient conditions for random asynchronous consensus as follows.

\begin{theo}\label{them:1}
The asynchronous DCA (\ref{eq:asynchronous}) generates consensus almost surely if all the following conditions hold:
\begin{itemize}
\item[a)] $\mathcal{G}(A)$ is rooted and $\sigma_k\in 2^V$.
\item[b)] There exists $\alpha>0$ such that if $\mathbb{P}(\sigma_{k}\,|\,\sigma_{(k-1):1})\neq 0$, then $\mathbb{P}(\sigma_{k}\,|\,\sigma_{(k-1):1})\geq \alpha$.
\item[c)] For any given past values of $\sigma_{(k-1):1}$, the set
    $$\mathscr{I}_{\sigma_{(k-1):1}}= \{\sigma: \mathbb{P}(\sigma_{k}=\sigma\,|\,\sigma_{(k-1):1})\neq 0\}$$
    only depends on $k$ but not the historic values $\sigma_{(k-1):1}$, i.e, there exists $\mathscr{I}_{k}$ such that
    $$\mathscr{I}_{\sigma_{(k-1):1}}=\mathscr{I}_{k},\quad \forall\, \sigma_{(k-1):1}\in \underbrace{2^V\times 2^V\cdots \times 2^V}_{k-1\,\, \mbox{\footnotesize times}}. $$
\item[d)] There exists $q>0$ such that
\begin{eqnarray}
\bigcup_{\tau=k}^{k+q-1}\left(\bigcup_{\sigma\in\mathscr{I}_\tau}\sigma\right)= V, \quad\forall\, k\geq 1. \nonumber
\end{eqnarray}
\item[e)] There exists a strongly connected component $\chi$ of $\mathbbm{r}(A)$ such that for any $j\in \chi$, it holds $\mathscr{I}_{k}^{j}\neq\emptyset$ and
    \begin{eqnarray}
    \chi\bigcap\left(\bigcap_{\sigma\in \mathscr{I}_{k}^{j} } \sigma\right) = j, \quad \forall\, k\geq 1,  \nonumber
    \end{eqnarray}
    where
    $$\mathscr{I}_{k}^{j}= \{\sigma: \sigma \in \mathscr{I}_k \mbox{ and } j\in \sigma\}.$$
\end{itemize}
\end{theo}

The meaning of the five conditions in {\em Theorem} \ref{them:1} is intuitive: condition a) gives the topological condition on $\mathcal{G}(A)$; condition b) requires a positive infimum on all the nonzero transition probabilities; condition c) is called historic independence of nonzero probabilities, which requires the set $\mathscr{I}_{\sigma_{(k-1):1}}$ to be independent of the historic values $\sigma_{(k-1):1}$; condition d) is called the joint coverage condition, which means each node of $V$ has a nonzero probability to be chosen to update in every consecutive $q$ steps; condition e) is called the quasi-singleton property of nodes in $\chi$, which means once some node $j\in \chi$ has a nonzero probability to be chosen at time $k$, then the intersection of $\chi$ and all the possible values of $\sigma_k$ which contain $j$ is $j$ itself.

To illustrate the meanings of the conditions c)-e) in {\em Theorem} \ref{them:1}, we give the following example:

\begin{exam}\label{exam:ExampleForTheorem1}
Suppose $\{\sigma_k\}_{k=1}^{\infty}$ are generated by i.i.d random variables, $V=\{1,2,3,4\}$, and
\begin{eqnarray}
A = \left(\begin{array}{cccc}
0& 1& 0& 0\\
0& 0& 1& 0\\
1& 0& 0& 0\\
0& 0& 1& 0
\end{array}\right), \label{eq:eg4}
\end{eqnarray}
where $\mathbbm{r}(A)=\{1,2,3\}$. Consider a subset of $2^{V}$, such as $$\mathscr{I}_k=\{\{1,2,4\},\{1,3,4\},\{2,3\}\}\footnote{The set $\mathscr{I}_k$ can be set as time-varying, but we only give a time-invariant example here for simplicity.}, \forall k\geq 1.$$
Since $\mathscr{I}_k$ does not rely on the historic values $\sigma_{(k-1):1}$, condition c) naturally holds. Note that
\begin{eqnarray}
& \mathscr{I}^{1}_{k} =  \{\{1,2,4\},\{1,3,4\}\}, \quad
\mathscr{I}^{2}_{k} =  \{\{2,3\},\{1,2,4\}\}, &\nonumber\\
& \mathscr{I}^{3}_{k} =  \{\{1,3,4\},\{2,3\}\}. & \nonumber
\end{eqnarray}
One can verify that
$$\bigcap_{\sigma\in\mathscr{I}^{1}_{k}}\sigma= \{1,4\}, \quad
\bigcap_{\sigma\in\mathscr{I}^{2}_{k}}\sigma= 2, \quad
\bigcap_{\sigma\in\mathscr{I}^{3}_{k}}\sigma= 3.$$
Since $\mathbbm{r}(A)=\{1,2,3\}$ and letting $\chi=\mathbbm{r}(A)$, condition e) of {\em Theorem} \ref{them:1} is satisfied. One further sets $q=1$ and calculates that $\bigcup_{\sigma\in {\mathscr{I}}_{k}}\sigma = \{1,2,3\} = V$.
Hence, condition d) holds too.
\end{exam}

The structure of the transition probability in conditions c)-e) of {\em Theorem} \ref{them:1} can be described by a trellis graph \cite{Automatica-Touri-Nedic}. The trellis graph of the random process $\{\Xi_k\}_{k\geq 1}$ is an infinite directed graph $\mathscr{T}=(\mathcal{V}, \mathcal{E}, \{\sigma_k\}_{k\geq 1})$, where $\mathcal{V}$ is the infinite grid $2^V\times \mathbb{Z}^+$ and
\begin{eqnarray}
\mathcal{E}= \{((\sigma,k),(\sigma',k+1))\,|\,\sigma,\sigma'\in 2^V,k\geq 1\}. \nonumber
\end{eqnarray}
According to the above definition, the link in a trellis graph is pointed from time $k$ to time $k+1$ if the corresponding transition probability is nonzero.



\begin{figure}[htbp]
\begin{center}
\resizebox{8.5cm}{2.5cm}{\includegraphics{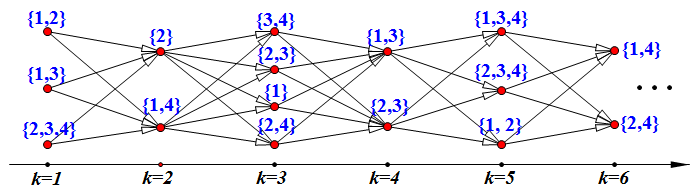}}
\end{center}
\caption{An example of trellis graph which satisfies the conditions of {\em Theorem} \ref{them:1}. }\label{fig:trellis}
\end{figure}

As shown in Fig. \ref{fig:trellis}, the given trellis graph satisfies the conditions of {\em Theorem} \ref{them:1} for $V=\{1,2,3,4\}$, $\mathbbm{r}(A)=\chi=\{1,2\}$, and $q=2$. It should be noted that the weight of each edge in Fig. \ref{fig:trellis} is dependent on the historic values.

In {\em Theorem} \ref{them:1}, we do not require $\{\Xi_k\}_{k\geq 1}$ to be i.i.d or Markovian. Hence, some popular methods for linear stochastic systems, such as LMI (linear matrix inequality), cannot be be easily applied. Even if $\{\Xi_k\}_{k\geq 1}$ are Markovian, the number of states in the Markovian chain is $|\Omega|=2^N$ and hence it generates $2^N$ LMIs by using the method given in \cite{MCSS-Boukas}. Unfortunately, LMIs with such a huge dimension is very difficult to be solved.


The conditions of {\em Theorem} \ref{them:1} imply the sample space of $A_{\sigma_k}$ may contain matrices with zero diagonal entries. Traditionally, the positivity of the diagonal entries in the stochastic matrices plays a very important role in random or deterministic consensus, which can be seen in \cite{I3ESAC-Shi}, \cite{I3EAC-Porfiri}, and \cite{I3EAC-Alireza}. In {\em Theorem} \ref{them:1}, we do not require positivity of the diagonal entries in the sample space $\mathfrak{A}$ to guarantee consensus of asynchronous DCA (\ref{eq:asynchronous}).

As simple implications of {\em Theorem} \ref{them:1}, we consider the following two practical cases of the random process $\{\Xi_k\}_{k\geq 1}$ which governs the asynchronous iteration:
\begin{itemize}
\item[a)] A global clock: In this case, a global clock determines which node to update, and there is only one node updates at each time. Hence, $\sigma_k\in V$ is the unique node updates at time $k$.
\item[b)] Independent asynchronous clocks: In this case, each node has a local clock, such a clock determines the update of the corresponding node independently. Hence, the updated nodes $\sigma_k$ at time $k$ is a set and it can be decomposed as $\sigma_k=\bigcup_{j=1}^N\theta_k^{(j)}$, where the $j$th random variable $\theta_k^{(j)}\in\{\{j\}, \emptyset\}$ is associated with the $j$th local clock. If $\theta_k^{(j)}=\{j\}$, then node $j$ updates at time $k$. If $\theta_k^{(j)}=\emptyset$, then node $j$ does not update at time $k$.
\end{itemize}

The sufficient conditions for almost sure consensus of asynchronous DCA (\ref{eq:asynchronous}) in the above two cases are given as follows.

\begin{theo}\label{them:globalasynchronousclock}
The asynchronous DCA (\ref{eq:asynchronous}) generates consensus almost surely if both of the following conditions hold:
\begin{itemize}
\item[a)] $\sigma_k\in V$ and $\mathcal{G}(A)$ is rooted.
\item[b)] There exists $\alpha>0$ such that for any $k\geq 1$
$$\mathbb{P}(\sigma_k=j)\geq \alpha,\quad \forall\, j\in V.$$
\end{itemize}
\end{theo}

\begin{theo}\label{them:independentclocks}
The asynchronous DCA (\ref{eq:asynchronous}) generates consensus almost surely if all the following conditions hold:
\begin{itemize}
\item[a)] $\mathcal{G}(A)$ is rooted.
\item[b)] For any $k\geq 1$, $\sigma_k=\bigcup_{j=1}^N\theta_k^{(j)}$, where $\theta_k^{(j)}\in\{\{j\},\emptyset\}$.
\item[c)] There exists $\alpha\in(0,\frac{1}{2}]$ such that for any $k\geq 1$, $$\mathbb{P}(\theta_k^{(j)}=\{j\})\in[\alpha, 1-\alpha],\quad \forall\, j\in V.$$
\end{itemize}
\end{theo}

\section{Technical Skeleton}\label{se:4}

We present the technical skeleton of the proof of {\em Theorem} \ref{them:1} as follows and the details can be found in the appendix.

First, we prove the equivalence between random asynchronous consensus and the convergence of the ergodic coefficient of products of stochastic matrices (see, {\em Lemma} \ref{lem:uniformconverge} and \ref{lem:scramblingconvergence}).

\begin{lemm}\label{lem:uniformconverge}
If the asynchronous DCA (\ref{eq:asynchronous}) realizes consensus almost surely under the conditions of {\em Theorem} \ref{them:1}, then it realizes consensus uniformly with respect to the region of the initial value $\mathscr{A}=\{x\,:\,\Delta(x)\leq 1\}$.
\end{lemm}

\begin{lemm}\label{lem:scramblingconvergence}
The asynchronous DCA (\ref{eq:asynchronous}) realize consensus almost surely if and only if $\lim_{k\rightarrow\infty}\mathbb{P}(\lambda(A_{\sigma(k:1)})\geq\varepsilon)=0$ holds for any $\varepsilon>0$.
\end{lemm}

Second, we show that the transition probability among $\{\sigma_k\}_{k\geq 1}$ has a special property: the backward probabilities in the trellis graph have a positive infimum (see, {\em Lemma} \ref{lem:backprob}).

\begin{lemm}\label{lem:backprob}
If the conditions of {\em Theorem} \ref{them:1} are satisfied, then for any integer $T\geq k$, there exists $\gamma\in(0,1)$ such that
\begin{eqnarray}
\mathbb{P}(\sigma_{k-1}\,|\, \sigma_{T:k})\geq \gamma, \quad \forall \,T\geq k\geq 2, \nonumber
\end{eqnarray}
when $\mathbb{P}(\sigma_{k-1}\,|\, \sigma_{T:k})\neq 0$.
\end{lemm}

Third, based on the result of step 2, the asynchronous consensus problem will be mapped to the problem of random backward walk along a labelled directed cycle (see, {\em Lemma} \ref{lem:rate}, \ref{lem:backwalk}), and the existence of such a cycle is guaranteed by {\em Lemma} \ref{prop:graphcycle}. Furthermore, we show that the random product of asynchronous matrix sequence has a nonzero probability to be scrambling (see, {\em Lemma} \ref{lem:sharingprob}, \ref{lem:scramblingpositive}).

\begin{lemm}\label{lem:rate}
Given a sequence of column stochastic matrices $\{P_{k}\}_{k=1}^{\infty}$, if $P_{k}\geq W$ and $P_k\sim W$ for each $k\geq 1$, and $\mathcal{G}(W^T)$ is rooted with node $1$ as the unique root which contains a self-loop, then there exist $c_0>0$ and $\beta\in(0,1)$ such that $\|P_k\cdots P_2P_1-\mathbbm{e}_1\mathbbold{1}^T\|\leq c_0\beta^k$,
where $\mathbbm{e}_1$ is the elementary unit vector with the $1$st entry being $1$ and the other entries being $0$.
\end{lemm}

\begin{lemm}\label{lem:backwalk}
(Random walk along a labelled directed cycle)
Given a labelled directed cycle $\mathscr{C}$ with the set of nodes $V=\{1,2,\cdots l\}$ and the corresponding set of labels $V_l=\{{\rm label}(i):i\in V\}$, we define the following random walk $\mathcal{W}$ along this cycle: Given two nodes moving along the cycle $\mathscr{C}$, the positions of them at time $k$ are denoted by $i_k$ and $j_k$ with the corresponding historic values $W_k=(i_\tau,j_\tau)_{\tau=1}^k$, and the transition from $(i_k,j_k)$ to $(i_{k+1},j_{k+1})$ are described by:

When ${\rm label}(i_k)\neq {\rm label}(j_k)$, it holds that
\begin{eqnarray}
\begin{array}{lclcl}
\mathbb{P}((i_{k+1},j_{k+1})&=&(i_k,j')\,|\,W_k) &\geq& \gamma, \\
\mathbb{P}((i_{k+1},j_{k+1})&=&(i',j_k)\,|\,W_k) &\geq& \gamma, \\
\mathbb{P}((i_{k+1},j_{k+1})&\in&\{(i_k,j_k),(i',j')\}\,|\,W_k) &\geq& \gamma,
\end{array} \nonumber
\end{eqnarray}
where $i'\rightarrow i_k$ and $j'\rightarrow j_k$ in $\mathscr{C}$, and
\begin{eqnarray}
\sum_{i\in \{i_k,i'\}, j\in\{j_k,j'\}}\mathbb{P}((i_{k+1},j_{k+1})=(i,j)\,|\,W_k)=1; \nonumber
\end{eqnarray}

When ${\rm label}(i_k)= {\rm label}(j_k)$, it holds
\begin{eqnarray}
\mathbb{P}((i_{k+1},j_{k+1})=(i_k,j_k)\,|\,W_k) =1.  \nonumber
\end{eqnarray}

For the above random walk $\mathcal{W}$, there exists an integer $k^*>0$ and real numbers $\mu_k\in(0,1)$ such that
\begin{eqnarray}
\mathbb{P}({\rm label}(i_{k+1})={\rm label}(j_{k+1})\,|\,W_k)\geq \mu_k. \nonumber
\end{eqnarray}
holds for any $k\geq k^*$.
\end{lemm}

\begin{lemm}\label{lem:sharingprob}
If the conditions of {\em Theorem} \ref{them:1} are satisfied, then given any two nodes $i,j\in V$ ($i\neq j$), there exists $T^*$ such that $\mathbb{P}(\mathcal{N}(A_{\sigma_{T:1}}, i)\cap \mathcal{N}(A_{\sigma_{T:1}}, j)\neq \emptyset) \,\geq\, h_T \, \in \,(0,1)$
holds for any $T\geq T^*$.
\end{lemm}

\begin{lemm}\label{lem:scramblingpositive}
If the conditions of {\em Theorem} \ref{them:1} are satisfied, then $\mathbb{P}(\lambda(A_{\sigma_{M:1}})\leq 1-\delta^M) \geq h_T^{\frac{N(N-1)}{2}}$, where $M=\frac{(N-1)N}{2}T$, $T\geq T^*$, $\delta$ is the minimal positive entry of $A$, $T^*$ and $h_T$ are given in {\em Lemma} \ref{lem:sharingprob}.
\end{lemm}

Finally, the proof of {\em Theorem} \ref{them:1} can be obtained with the assistance of step 1 by partitioning the matrix product $A_{M:1}$ to many subproducts with each having a nonzero probability to be scrambling. By using the basic properties of ergodic  coefficient, the convergence can be derived if scrambling matrices appear for infinitely many times.

The proof of the above lemmas and the proof of {\em Theorem} \ref{them:1} are given in the appendix.
In particular, both {\em Theorem} \ref{them:globalasynchronousclock} and \ref{them:independentclocks} can be directly obtained from {\em Theorem} \ref{them:1}, hence the proofs of them have been omitted. Moreover, one can refer to the appendix for more corollaries of {\em Theorem} \ref{them:1}.


\section{Numerical Examples}\label{se:5}

In this section, we give two numerical examples to verify the effectiveness of {\em Theorem} \ref{them:globalasynchronousclock} and \ref{them:independentclocks}.

\begin{figure}[htbp]
\begin{center}
\resizebox{7.0cm}{5.0cm}{\includegraphics{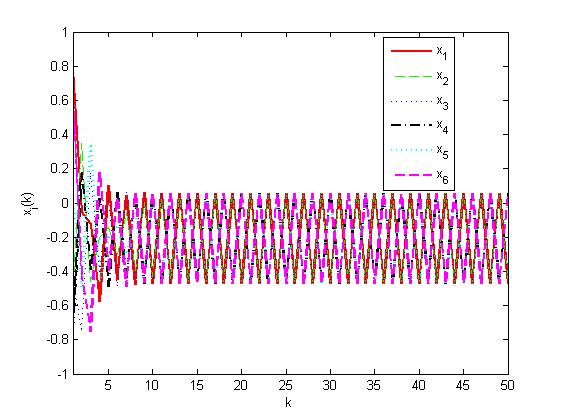}}
\end{center}
\caption{Synchronous DCA (\ref{eq:synchronous}) with $A$ given in (\ref{eq:A}): cannot realize consensus. }\label{fig:0}
\end{figure}

\begin{figure}[htbp]
\begin{center}
\resizebox{7.0cm}{5.0cm}{\includegraphics{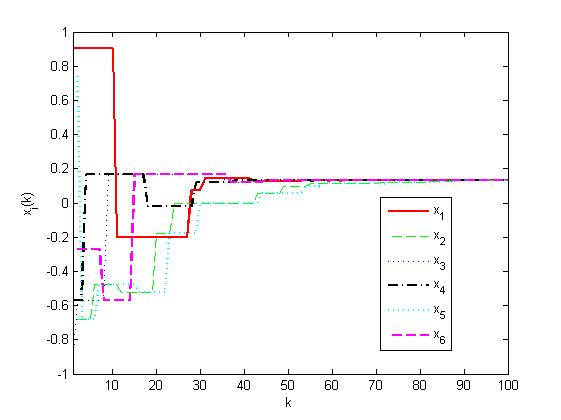}}
\end{center}
\caption{Asynchronous DCA (\ref{eq:asynchronous}) with $A$ given in (\ref{eq:A}): a global clock. }\label{fig:1}
\end{figure}

Consider the following stochastic matrix
\begin{eqnarray}
A = \left(\begin{array}{cccccc}
0 & 0 & 0 & 0.5 & 0 & 0.5\\
0 & 0 & 0 & 0 & 0.5 & 0.5\\
0 & 0 & 0 & 1 & 0 & 0\\
0.5 & 0 & 0.5 & 0 & 0 & 0\\
0 & 1 & 0 & 0 & 0 & 0\\
0 & 0 & 1 & 0 & 0 & 0
\end{array}\right), \label{eq:A}
\end{eqnarray}
where $\mathcal{G}(A)$ is strongly connected. By using the method given in \cite{SCL-Chen}, one can see $A$ is not an SIA matrix. As shown in Fig. \ref{fig:0}, the synchronous DCA (\ref{eq:synchronous}) cannot realize consensus with the above $A$. However, when we implement the asynchronous DCA (\ref{eq:asynchronous}) with the same matrix $A$, the dynamics becomes quite different:

Firstly, suppose DCA (\ref{eq:asynchronous}) has a global clock. For any $j\in\{1,2,\cdots, 6\}$, one sets $\mathbb{P}(\sigma_k=j)=1/6$ and the conditions of {\em Theorem} \ref{them:globalasynchronousclock} (also, {\em Theorem} \ref{them:1}) are satisfied. As shown in Fig. \ref{fig:1}, the agents realize consensus almost surely when the initial values are chosen randomly from $[-1,1]$.

\begin{figure}[htbp]
\begin{center}
\resizebox{7.0cm}{5.0cm}{\includegraphics{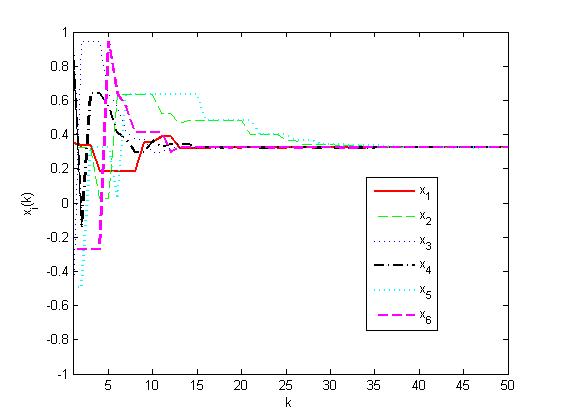}}
\end{center}
\caption{Asynchronous DCA (\ref{eq:asynchronous}) with $A$ given in (\ref{eq:A}): independent asynchronous clocks. }\label{fig:2}
\end{figure}

Next, suppose DCA (\ref{eq:asynchronous}) has independent local clocks. For any $j\in\{1,2,\cdots, 6\}$, one sets $\mathbb{P}(\theta_k^{(j)}=j)=\mathbb{P}(\theta_k^{(j)}=\phi)=\frac{1}{2}$ and the conditions of {\em Theorem} \ref{them:independentclocks} (also, {\em Theorem} \ref{them:1}) are satisfied. As shown in Fig. \ref{fig:2}, the agents realize consensus almost surely when the initial values are chosen randomly from $[-1,1]$.

\section{Conclusions}\label{se:6}

This paper has investigated the asynchronous iteration problem of distributed coordination algorithms in the stochastic setting.  We have found that the topology only needs to be rooted for consensus of DCA with asynchronous iteration in the given stochastic setting, which is in sharp contrast to the deterministic setting where the topology should be rooted and aperiodic. In the future, we will apply the proposed DCA with random asynchronous iteration to resolve practical engineering problems, such as clock synchronization in wireless sensor networks.





\appendix\label{appendix}

\subsection{Examples}\label{appendix:examples}

\begin{exam} \label{exam:NonConvergentAsynchronousSeqWithSIAMatrix}
Given
\begin{eqnarray}
A = \left(\begin{array}{ccccc}
0& 0& 0& 0& 1\\
1& 0& 0& 0& 0\\
0.5& 0.5& 0& 0& 0\\
0& 0& 1& 0& 0\\
0& 0& 0& 1& 0
\end{array}\right), \nonumber
\end{eqnarray}
which is SIA, but the corresponding asynchronous iteration sequence $\{A_5, A_4, A_1, A_2, A_3, A_5, A_4, A_1, A_2, A_3, \cdots\}$ does not generate consensus since
\begin{eqnarray}
A_3A_2A_1A_4A_5 = \left(\begin{array}{ccccc}
0& 0& 0& 1& 0\\
0& 0& 0& 1& 0\\
0& 0& 0& 1& 0\\
0& 0& 1& 0& 0\\
0& 0& 0& 1& 0
\end{array}\right) \nonumber
\end{eqnarray}
is not SIA.
\end{exam}

\begin{exam}\label{exam:ConvergentAsynchronousSeqWithNonSIAMatrix}
Given
\begin{eqnarray}
A = \left(\begin{array}{cc}
0& 1\\
1& 0
\end{array}\right), \nonumber
\end{eqnarray}
which is not SIA but the corresponding asynchronous iteration sequence $\{A_2,A_1,A_2,A_1,\cdots\}$ gives rise to consensus since
\begin{eqnarray}
A_1A_2 = \left(\begin{array}{cc}
1& 0\\
1& 0
\end{array}\right) \nonumber
\end{eqnarray}
is SIA.
\end{exam}

\subsection{Technical Details}
The proofs of the lemmas are given as follows.

{ \textbf{Proof of Lemma \ref{lem:uniformconverge}}: } Without loss of generality, suppose $x=(x_1,x_2,\cdots,$ $ x_N)^T\in \mathscr{A}$ and $\underline{x}=\min_{i=1}^N x_i$, then $\tilde{x}=x-\underline{x}$ also satisfies $\tilde{x}\in\mathcal{A}$. Particularly, $\min_{i=1}^N \tilde{x}_i=0$ and $\max_{i=1}^N \tilde{x}_i\leq 1$, where $\tilde{x}_i=x_i-\underline{x}$. Hence, one can write $x$ in the form of
\begin{eqnarray}
x = \underline{x}\mathbbold{1} + \sum_{i=1}^N\lambda_i \mathbbm{e}_i. \nonumber
\end{eqnarray}
where $0\leq \lambda_i\leq 1$ for any $1\leq i\leq N$, $\mathbbm{e}_i$ is the $i$th elementary vector.

Consider a set of initial values $\{\mathbbm{e}_i\}_{i=1}^N$. Then each trajectory starting from $\mathbbm{e}_i$ satisfies (\ref{eq:equalform}). For any given $\varepsilon>0$ and $\eta>0$, there exists $N_i$ such that
\begin{eqnarray}
\mathbb{P}\left(\|Px(k)\|\geq \frac{\varepsilon}{N}\,|\, x(1)=\mathbbm{e}_i\right)\leq \eta
\end{eqnarray}
for any $k\geq N_i$. Define $\widetilde{N}=\max_{i=1}^N N_i$. Then
\begin{eqnarray}
\mathbb{P}\left(\|PA_{\sigma_{k:1}}\mathbbm{e}_i\|\geq \frac{\varepsilon}{N}\right)\leq \eta
\end{eqnarray}
for any $k\geq \widetilde{N}$ and $i\in V$. For any initial value $x(1)=x\in \mathscr{A}$, one has
\begin{eqnarray}
\|Px(k+1)\| &=& \|PA_{\sigma_{k:1}}x(1)\| \nonumber\\
&\leq& \sum_{i=1}^N\lambda_i \|PA_{\sigma_{k:1}}\mathbbm{e}_i\| \nonumber\\
&\leq& \sum_{i=1}^N \|PA_{\sigma_{k:1}}\mathbbm{e}_i\| \nonumber
\end{eqnarray}
and
\begin{eqnarray}
\mathbb{P}(\|Px(k+1)\|\geq \varepsilon) &\leq & \mathbb{P}\left(\sum_{i=1}^N \|PA_{\sigma_{k:1}}\mathbbm{e}_i\|\geq \varepsilon\right) \nonumber\\
&\leq & \max_{i=1}^N  \mathbb{P}\left( \|PA_{\sigma_{k:1}}\mathbbm{e}_i\|\geq \frac{\varepsilon}{N}\right). \nonumber
\end{eqnarray}
Note that $\mathbb{P}\left( \|PA_{\sigma_{k:1}}\mathbbm{e}_i\|\geq \frac{\varepsilon}{N}\right)\leq \eta$ for each $i\in V$ when $k\geq \widetilde{N}$, and the proof is hence completed. \hfill $\square$

\textbf{Proof of Lemma \ref{lem:scramblingconvergence}: } On one hand, for any $x=(x_1,x_2,\cdots, x_N)^T\in\mathbb{R}^N$, there is
\begin{eqnarray}
\|Px\| &=& \sqrt{\sum_{j=1}^N\left(x_j-\frac{1}{N}\sum_{i=1}^N x_i\right)^2}\nonumber\\
&=& \sqrt{\frac{1}{N^2}\sum_{j=1}^N\left(\sum_{i=1}^N (x_j-x_i)\right)^2}\nonumber\\
&\leq & \sqrt{\frac{1}{N^2}\sum_{j=1}^N\left(\sum_{i=1}^N |x_j-x_i|\right)^2} \nonumber\\
&\leq & \sqrt{N}\cdot\Delta(x), \nonumber
\end{eqnarray}

On the other hand, suppose that $\overline{x}=\max_{i=1}^N x_i$ and $\underline{x}=\min_{i=1}^N x_i$, there is
\begin{eqnarray}
\|Px\| &=& \sqrt{\sum_{j=1}^N\left(x_j-\frac{1}{N}\sum_{i=1}^N x_i\right)^2}\nonumber\\
&\geq & \sqrt{\left(\overline{x}-\frac{1}{N}\sum_{i=1}^N x_i\right)^2 + \left(\underline{x}-\frac{1}{N}\sum_{i=1}^N x_i\right)^2} \nonumber\\
&\geq& \frac{1}{\sqrt{2}}\cdot \Delta(x). \nonumber
\end{eqnarray}

Hence,
\begin{eqnarray}
\frac{1}{\sqrt{2}}\cdot \Delta(x) \leq \|Px\|\leq \sqrt{N}\cdot \Delta(x). \label{eq:normeq}
\end{eqnarray}

Without loss of generality, suppose that $\Delta(x(1))=1$. We examine separately sufficiency and necessity as follows:

\emph{Sufficiency}: Since
$$\Delta(x(k+1))\leq \lambda(A_{\sigma_{k:1}})\Delta(x(1))= \lambda(A_{\sigma_{k:1}}),$$
one knows
\begin{eqnarray}
\mathbb{P}(\Delta(x(k+1))\geq \varepsilon) \leq \mathbb{P}(\lambda(A_{\sigma_{k:1}})\geq \varepsilon).  \nonumber
\end{eqnarray}
According to (\ref{eq:normeq}), one further knows that
\begin{eqnarray}
\mathbb{P}(\|Px(k+1)\|\geq \sqrt{N}\varepsilon) \leq \mathbb{P}(\Delta(x(k+1))\geq \varepsilon). \nonumber
\end{eqnarray}
Hence, $\lim_{k\rightarrow\infty}\mathbb{P}(\lambda(A_{\sigma_{k:1}})\geq \varepsilon)=0$ implies $\lim_{k\rightarrow\infty}\mathbb{P}$ $(\|Px(k+1)\|\geq \sqrt{N}\varepsilon)=0$. Then, sufficiency follows from the fact that $\epsilon$ is chosen arbitrarily small.

\emph{Necessity}: According to (\ref{eq:normeq}), one knows that
$$\mathbb{P}(\Delta(x(k+1))\geq \sqrt{2}\varepsilon)\leq \mathbb{P}(\|Px(k+1)\|\geq\varepsilon).$$
Hence, $\lim_{k\rightarrow\infty}\mathbb{P}(\|Px(k+1)\|\geq\varepsilon)=0$ implies $\lim_{k\rightarrow\infty}\mathbb{P}$ $(\Delta(A_{\sigma_{k:1}}x(1))\geq \sqrt{2}\varepsilon)=0$. Based on {\em Lemma} \ref{lem:uniformconverge} and the fact that
$$\lambda(A_{\sigma_{k:1}}) = \sup_{\Delta(x(1))=1}\Delta(A_{\sigma_{k:1}}x(1)),$$
one knows $\lim_{k\rightarrow\infty}\mathbb{P}(\lambda(A_{\sigma_{k:1}})\geq \sqrt{2}\,\varepsilon)=0$. Necessity also follows from the arbitrary choice of $\epsilon$. \hfill $\square$

\textbf{Proof of Lemma \ref{lem:backprob}: } According to condition b) of {\em Theorem} \ref{them:1}, for any $k\geq 1$ and $j\in V$, there is
\begin{eqnarray}
\mathbb{P}(\sigma_k)
&=& \sum_{\sigma_{1},\cdots, \sigma_{k-1}}\mathbb{P}(\sigma_{k}\,|\, \sigma_{(k-1):1})\mathbb{P}(\sigma_{(k-1):1})  \nonumber\\
&\geq& \alpha\cdot \sum_{\mathbb{P}(\sigma_{k}\,|\, \sigma_{(k-1):1})\neq 0}\mathbb{P}(\sigma_{(k-1):1}). \nonumber
\end{eqnarray}
Note that condition c) of {\em Theorem} \ref{them:1} implies that
\begin{eqnarray}
& & \{\sigma_{(k-1):1}\,|\, \mathbb{P}(\sigma_{k}=\sigma\,|\, \sigma_{(k-1):1})\neq 0\} \nonumber\\
&=& \{\sigma_{(k-1):1}\,|\, \mathbb{P}(\sigma_{(k-1):1})\neq 0\} \nonumber
\end{eqnarray}
for any fixed $\sigma\in\mathscr{I}_k$, based on which one derives
\begin{eqnarray}
\mathbb{P}(\sigma_k)
&\geq& \alpha\cdot \sum_{\mathbb{P}(\sigma_{(k-1):1})\neq 0}\mathbb{P}(\sigma_{(k-1):1}) =\alpha.  \nonumber
\end{eqnarray}

Similarly, one obtains
\begin{eqnarray}
& & \mathbb{P}(\sigma_{(k+1):k}) \nonumber\\
&=& \sum_{\sigma_{1},\cdots, \sigma_{k-1}}\mathbb{P}(\sigma_{k+1}\,|\, \sigma_{k:1})\mathbb{P}(\sigma_{k}\,|\, \sigma_{(k-1):1})\mathbb{P}(\sigma_{(k-1):1}) \nonumber\\
&\geq& \alpha^2. \nonumber
\end{eqnarray}
Repeating the above procedure one obtains that if $\mathbb{P}(\sigma_{T:k})\neq 0$, then
\begin{eqnarray}
\mathbb{P}(\sigma_{T:k}) \geq \alpha^{T-k+1} \nonumber
\end{eqnarray}
for any $T\geq k\geq 1$.

Based on the above discussion, when $\mathbb{P}(\sigma_{k-1}\,|\, \sigma_{T:k})\neq 0$, it follows
\begin{eqnarray}
\mathbb{P}(\sigma_{k-1}\,|\, \sigma_{T:k})
= \frac{\mathbb{P}(\sigma_{T:(k-1)})}{\mathbb{P}(\sigma_{T:k})}
\geq \alpha^{T-k+2}
\geq \alpha^{T} = \gamma, \nonumber
\end{eqnarray}
and the proof is hence completed. \hfill $\square$

\textbf{Proof of Lemma \ref{lem:rate}: } By using the result on convergence rate of consensus in \cite{TNN-Chen}, one can easily derive this result. \hfill $\square$

\begin{rema}
Lemma \ref{lem:backwalk} is critical for the development of this paper. It describes the following phenomenon: given a labelled directed cycle $\mathscr{C}$ and two nodes in this cycle, the movement of these two nodes is governed by a so-called backward random walk $\mathcal{W}$; {\em Lemma} \ref{lem:backwalk} shows that these two nodes have a nonzero probability to share the same label after several steps of iteration. The proof of {\em Lemma} \ref{lem:backwalk} is based on this comparison technique: we construct another random walk $\mathcal{W}'$ and show that the two nodes have a higher probability to share the same label under the effect of $\mathcal{W}$ than to collide under the effect of $\mathcal{W}'$ (see (\ref{eq:cycle1})); then we show that $\mathcal{W}'$ makes the initial two nodes collide along $\mathscr{C}$ after several iterations (see (\ref{eq:cycle3})); we finally finish the proof by combining these two facts.
\end{rema}

\textbf{Proof of Lemma \ref{lem:backwalk}: } Define the set
\begin{eqnarray}
\mathcal{T} &=& \left\{(g_k,h_k)_{k=1}^m\,|\, g_k,h_k\in V_{l},\, {\rm label}(g_k)\neq {\rm label}(h_k) \right. \nonumber\\
& & \left. \qquad \qquad \qquad \mbox{holds for any } 1\leq k\leq m, \,m\geq 1\right\}. \nonumber
\end{eqnarray}
Using this set and the random walk $\mathcal{W}$, we construct another random walk $\mathcal{W}'$ along $\mathscr{C}$: at each time $k$, the positions of two nodes are denoted by $i_k'$ and $j_k'$, and the corresponding historic values are denoted by $W_k'=(i_\tau',j_\tau')_{\tau=1}^k$. The transition from time $k$ to $k+1$ follows:

When $i_k'\neq j_k'$ and $W_k'\in \mathcal{T}$, it holds
\begin{eqnarray}
\mathbb{P}((i_{k+1}',j_{k+1}')=(i_{k+1},j_{k+1})\,|\,W_k')=1. \nonumber
\end{eqnarray}

When $i_k'\neq j_k'$ and $W_k'\notin \mathcal{T}$, it follows
\begin{eqnarray}
\begin{array}{lclcl}
\mathbb{P}((i_{k+1}',j_{k+1}')&=&(i_k',j')\,|\,W_k') &\geq& \gamma, \\
\mathbb{P}((i_{k+1}',j_{k+1}')&=&(i',j_k')\,|\,W_k') &\geq& \gamma, \\
\mathbb{P}((i_{k+1}',j_{k+1}')&\in&\{(i_k',j_k'),(i',j')\}\,|\,W_k') &\geq& \gamma,
\end{array} \nonumber
\end{eqnarray}
where $i'\rightarrow i_k'$ and $j'\rightarrow j_k'$ in $\mathscr{C}$, and
\begin{eqnarray}
\sum_{i\in \{i_k',i'\}, j\in\{j_k',j'\}}\mathbb{P}((i_{k+1}',j_{k+1}')=(i,j)\,|\,W_k')=1; \nonumber
\end{eqnarray}

When $i_k'=j_k'$, we have
\begin{eqnarray}
\mathbb{P}((i_{k+1}',j_{k+1}')=(i_k',j_k')\,|\,W_k') = 1. \nonumber
\end{eqnarray}

We would like to point out that for any two nodes $i_1,j_1\in V_l$, it must be true that
\begin{eqnarray}
& &\mathbb{P}({\rm label}(i_{k+1})={\rm label}(j_{k+1})\,|\,i_1, j_1,\mathcal{W})\nonumber\\
&\geq &
\mathbb{P}(i_{k+1}'=j_{k+1}'\,|\,i_1, j_1,\mathcal{W}'). \label{eq:cycle1}
\end{eqnarray}
In fact, we can introduce
\begin{eqnarray}
\mathcal{F}_{k+1}&=& \left\{(i_\tau,j_\tau)_{\tau=1}^{k+1}\,|\, {\rm label}(i_{k+1})={\rm label}(j_{k+1}), \right.\nonumber\\
& & \left. (i_k,j_k) \mbox{ are generated by } \mathcal{W}\right\}, \nonumber\\
\mathcal{F}_{k+1}'&=& \left\{(i_\tau,j_\tau)_{\tau=1}^{k+1}\,|\, i_{k+1}=j_{k+1}, \right.\nonumber\\
& & \left. (i_k,j_k) \mbox{ are generated by } \mathcal{W}'\right\}, \nonumber
\end{eqnarray}
and construct a map $\mathscr{F}: \mathcal{F}_{k+1}'\rightarrow\mathcal{F}_{k+1}$ as follows: consider any $W_{k+1}'\in\mathcal{F}_{k+1}'$ and suppose $s\in[1,k+1]$ is the first time which satisfies ${\rm label}(i_s)={\rm label}(j_s)$; then define
\begin{eqnarray}
\mathscr{F}(W_{k+1}') &=& \{(i_1,j_1),\cdots (i_{s-1}, j_{s-1}), \underbrace{(i_s,j_s), \cdots (i_s,j_s)}_{k+2-s \mbox{\footnotesize\, times}}\} \nonumber\\
&\in& \mathcal{F}_{k+1}, \nonumber
\end{eqnarray}
which obviously indicates $\mathscr{F}(\mathcal{F}_{k+1}')\subseteq \mathcal{F}_{k+1}$.

For $w'\in \mathscr{F}(\mathcal{F}_{k+1}')$, the set $\mathscr{F}^{-1}(w')$ is a subset of $\mathcal{F}_{k+1}'$. Define $\mathscr{F}^{-1}(w_\zeta)$ ($\zeta\in\Upsilon, w_\zeta\in\mathcal{F}_{k+1}$) as the different sets with the form of $\mathscr{F}^{-1}(w')$ ($w'\in \mathscr{F}(\mathcal{F}_{k+1}')$). Then
\begin{eqnarray}
\bigcup_{\zeta\in\Upsilon}\mathscr{F}^{-1}(w_{\zeta}) = \mathcal{F}_{k+1}' \nonumber
\end{eqnarray}
and $\mathscr{F}^{-1}(w_{\zeta_1})\bigcap \mathscr{F}^{-1}(w_{\zeta_2})=\emptyset$ for $\zeta_1\neq\zeta_2$. According to the definitions of $\mathcal{W}$ and $\mathcal{W}'$, it holds
\begin{eqnarray}
\sum_{w'\in \mathscr{F}^{-1}(w_\zeta)}\mathbb{P}(w'\,|\,i_1,j_1,\mathcal{W}) \leq \mathbb{P}(w_{\zeta}\,|\,i_1,j_1,\mathcal{W}),\quad \forall\,\zeta\in\Upsilon.  \nonumber
\end{eqnarray}
Summarizing the above both sides w.r.t $\Upsilon$, one obtains that
\begin{eqnarray}
\sum_{w'\in \mathcal{F}_{k+1}'}\mathbb{P}(w'\,|\,i_1,j_1,\mathcal{W}) \leq \sum_{\zeta\in\Upsilon}\mathbb{P}(w_{\zeta}\,|\,i_1,j_1,\mathcal{W}).\nonumber
\end{eqnarray}
Noting that
\begin{eqnarray}
\sum_{\zeta\in\Upsilon}\mathbb{P}(w_{\zeta}\,|\,i_1,j_1,\mathcal{W})\leq
\sum_{w\in \mathcal{F}_{k+1}}\mathbb{P}(w'\,|\,i_1,j_1,\mathcal{W}), \nonumber
\end{eqnarray}
one obtains the result of (\ref{eq:cycle1}).

Given any pair of nodes $i, j\in V_l$ ($i\neq j$) in $\mathscr{C}$, we define ${\rm d}_{\mathscr{C}}(i,j)$ as the distance from $i$ to $j$ along the cycle $\mathscr{C}$. For example, in Fig. \ref{fig:cycle}, ${\rm d}_{\mathscr{C}}(1,2)=1$ but ${\rm d}_{\mathscr{C}}(2,1)=5$. In detail, ${\rm d}_{\mathscr{C}}(\cdot,\cdot)$ has the following properties
\begin{itemize}
\item[a)] ${\rm d}_{\mathscr{C}}(i,i)=0, \quad \forall\, i\in V$;
\item[b)] $0\leq {\rm d}_{\mathscr{C}}(i,j)\leq l-1, \quad \forall\, i,j\in V$;
\item[c)] ${\rm d}_{\mathscr{C}}(i,j)+{\rm d}_{\mathscr{C}}(i,j)=l$ if $i\neq j$.
\end{itemize}

We further denote
$${\rm d}_k=\mathbbm{d}_{\mathscr{C}}(i_k',j_k'),\quad k\geq 1,$$
where $i_k'$ and $j_k'$ are generated by $\mathcal{W}'$ with the initial values $i_1'$ and $j_1'$.  Based on ${\rm d}_k$ we make the following discussions:

When ${\rm d}_k=0$, we have
$\mathbb{P}({\rm d}_{k+1}=0\,|\, {\rm d}_k=0)=1$;

When ${\rm d}_k\neq 0$ and $(i_{k+1},j_{k+1})\in\{(i_k,j_k),(i',j')\}$, we have
$\mathbb{P}({\rm d}_{k+1}={\rm d}_{k}\,|\, {\rm d}_k\neq 0)\geq \gamma$;

When $0<{\rm d}_k<l-1$ and $(i_{k+1},j_{k+1})=(i',j_k)$, we have
$\mathbb{P}({\rm d}_{k+1}={\rm d}_k+1\,|\, 0<{\rm d}_k<l-1)\geq \gamma$;

When ${\rm d}_k=l-1$ and $(i_{k+1},j_{k+1})=(i',j_k)$, we have
$\mathbb{P}({\rm d}_{k+1}=0\,|\, {\rm d}_k=l-1)\geq \gamma$;

When $1<{\rm d}_k<l$ and $(i_{k+1},j_{k+1})=(i_k,j')$, we have
$\mathbb{P}({\rm d}_{k+1}={\rm d}_k-1\,|\, 1<{\rm d}_k<l)\geq \gamma$;

When ${\rm d}_k=1$ and $(i_{k+1},j_{k+1})=(i_k,j')$, we have
$\mathbb{P}({\rm d}_{k+1}=0\,|\, {\rm d}_k=1)\geq \gamma$.

Summarizing the above discussions, one obtains that
\begin{eqnarray}
\mathbb{P}({\rm d}_{{k}+1}={\rm d}_{k}+1\,|\,0<{\rm d}_{k}<l-1)
&\geq& \gamma, \nonumber\\
\mathbb{P}({\rm d}_{{k}+1}=0\,|\, {\rm d}_{k}=l-1)
&\geq&  \gamma, \nonumber\\
\mathbb{P}({\rm d}_{{k}+1}={\rm d}_{k}-1\,|\,1<{\rm d}_{k}<l)
&\geq&  \gamma, \nonumber\\
\mathbb{P}({\rm d}_{{k}+1}=0\,|\, {\rm d}_{k}=1)
&\geq& \gamma, \nonumber\\
\mathbb{P}({\rm d}_{{k}+1}={\rm d}_{k}\,|\,{\rm d}_{k}\neq 0)
&\geq& \gamma, \nonumber\\
\mathbb{P}({\rm d}_{{k}+1}={\rm d}_{k}\,|\,{\rm d}_{k}= 0)
&=& 1, \nonumber
\end{eqnarray}
where the other conditional probabilities except the above cases are all zero.

Define
\begin{eqnarray}
\xi_\nu(k) &=& \mathbb{P}({\rm d}_{k}=\nu), \quad \forall\, \nu\in \{0,1,2,\cdots, l-1\}, \nonumber\\
\xi(k) &=& (\xi_0(k),\xi_1(k),\cdots, \xi_{l-1}(k))^{T}. \nonumber
\end{eqnarray}
The dynamics of $\xi(k)$ can be written as
\begin{eqnarray}
\xi(k+1)=P_k\xi(k), \nonumber
\end{eqnarray}
where $P_k\in\mathbb{R}^{l\times l}$ is a column stochastic matrix which satisfies $P_k\sim W$ and
\begin{eqnarray}
P_k \geq W = \left(\begin{array}{ccccccc}
1 & \gamma & 0 & 0 & \hdots & 0& \gamma\\
0 & \gamma & \gamma & 0 & \hdots&0 & 0\\
0 & \gamma & \gamma & \gamma & \hdots&0 & 0\\
0 & 0           & \gamma   & \gamma & \hdots&0 & 0\\
\vdots & \vdots  & \vdots   & \vdots & \ddots & \vdots & \vdots\\
0 & 0  & 0   & 0 &\hdots & \gamma&\gamma\\
0 & 0  & 0   & 0 &\hdots &\gamma &\gamma\\
\end{array}\right). \nonumber
\end{eqnarray}
Specifically, $\mathcal{G}(W^T)$ is rooted and node $1$ is the unique root. According to {\em Lemma} \ref{lem:rate}, one knows that $\lim_{k\rightarrow\infty}P_kP_{k-1}\cdots P_1=\mathbbm{e}_1\mathbbold{1}^T$ with an exponential rate. Note that for any fixed $i,j\in V_l$ ($i\neq j$), the vector $\xi(1)$ satisfies $\mathbbold{1}^T\xi(1)=1$. Since
\begin{eqnarray}
\xi(k+1)-\mathbbm{e}_1 = (P_k\cdots P_2P_1-\mathbbm{e}_1\mathbbold{1}^T)(\xi(1)-\mathbbm{e}_1), \nonumber
\end{eqnarray}
one derives the existence of $c_0>0$ and $\beta\in(0,1)$ such that
$|\xi_0(k)-1|\leq c_0\beta^k$ for any $k\geq 1$. Hence,
\begin{eqnarray}
\mathbb{P}({\rm d}_{k}=0) = \xi_0(k)\geq 1-|\xi_0(k)-1|\geq 1-c_0\beta^{k}, \nonumber
\end{eqnarray}
which implies that if $k \geq k^*=-\frac{\log c_0}{\log \beta}$, it holds
\begin{eqnarray}
\mathbb{P}({\rm d}_{k}=0) \geq 1-c_0\beta^{k} =\mu_k \in (0,1). \label{eq:cycle3}
\end{eqnarray}

Noting that
$$\mathbb{P}(i_{k+1}'=j_{k+1}'\,|\,i_1,j_1,\mathcal{W}')= \mathbb{P}({\rm d}_{k}=0\,|\,i_1,j_1,\mathcal{W}')$$
and applying (\ref{eq:cycle3}) to (\ref{eq:cycle1}), the proof is hence completed. \hfill $\square$

Based on the above lemma, one derives the following result.

\textbf{Proof of Lemma \ref{lem:sharingprob}: } Since $\mathcal{G}(A)$ is rooted and the graph induced by $\chi\subseteq\mathbbm{r}(A)$ is strongly connected, according to {\em Lemma} \ref{prop:graphcycle}, one can construct a labelled directed cycle $\mathscr{C}$ with length
$l\leq |\chi|(|\chi|-1)$ containing all the nodes of $\chi$. The set of nodes of $\mathscr{C}$ is $V=\{1,2,\cdots, l\}$ and the corresponding set of labels is $V_l=\{1,2,\cdots, |\chi|\}$.

Consider the subsequence $\{A_{\sigma_k}\}_{k=1}^T$ and any pair of nodes $i,j\in V$ $(i\neq j)$. According to {\em Lemma} \ref{lem:backprob} and condition d) of {\em Theorem} \ref{them:1}, one knows that
\begin{eqnarray}
\mathbb{P}(\chi\cap \mathcal{N}(A_{T:(T-(N-|\chi|)q+1)},i)\neq \emptyset \mbox{ and } \qquad\qquad\qquad \nonumber\\
\chi\cap \mathcal{N}(A_{T:(T-(N-|\chi|)q+1)},j)\neq \emptyset) \geq \gamma^{(N-|\chi|)q}.  \nonumber
\end{eqnarray}
Next, one defines
\begin{eqnarray}
\widetilde{T} = T-(N-|\chi|)q \nonumber
\end{eqnarray}
and chooses
\begin{eqnarray}
\hat{i} \in \chi\cap \mathcal{N}(A_{T:(\widetilde{T}+1)},i), \quad
\hat{j} \in \chi\cap \mathcal{N}(A_{T:(\widetilde{T}+1)},j). \nonumber
\end{eqnarray}
Based on $\hat{i}$ and $\hat{j}$ one chooses two nodes $i_{\widetilde{T}},j_{\widetilde{T}}\in \mathscr{C}$ with ${\rm label}(i_{\widetilde{T}})=\hat{i}$ and ${\rm label}(j_{\widetilde{T}})=\hat{j}$.

Based on the above pair of nodes $i_{\widetilde{T}}$, $j_{\widetilde{T}}$ and for any $k\leq \widetilde{T}$, we define iteratively
\begin{eqnarray}
i_{k-1} = \mathcal{N}_{\mathscr{C}}(i_k\,;\,\sigma_k) = \left\{\begin{array}{ll}
i_k, & \mbox{if } {\rm label}(i_k)\notin \sigma_k, \\
a', & \mbox{if } {\rm label}(i_k)\in \sigma_k \\
 & \,\,\,\mbox{ and } a'\rightarrow i_k \mbox{ in } \mathscr{C},
\end{array}\right. \nonumber\\
j_{k-1} = \mathcal{N}_{\mathscr{C}}(j_k\,;\,\sigma_k) = \left\{\begin{array}{ll}
j_k, & \mbox{if } {\rm label}(j_k)\notin \sigma_k, \\
a', & \mbox{if } {\rm label}(j_k)\in \sigma_k \\
 & \,\,\,\mbox{ and } a'\rightarrow j_k \mbox{ in } \mathscr{C}.
\end{array}\right. \nonumber
\end{eqnarray}
The compact form of the above iteration can be written as
\begin{eqnarray}
i_{k-1} &=& \mathcal{N}_{\mathscr{C}}(i_{\widetilde{T}}\,;\,\sigma_{{\widetilde{T}}:k}), \quad 1\leq k\leq {\widetilde{T}},\nonumber\\
j_{k-1} &=& \mathcal{N}_{\mathscr{C}}(j_{\widetilde{T}}\,;\,\sigma_{{\widetilde{T}}:k}), \quad 1\leq k\leq {\widetilde{T}}.\nonumber
\end{eqnarray}
As a simple illustration, one has
\begin{eqnarray}
6 = \mathcal{N}_{\mathscr{C}}(1\,;\,\sigma_k=\{1,3\}), \quad
1 = \mathcal{N}_{\mathscr{C}}(1\,;\,\sigma_k=\{2,4\}), \nonumber
\end{eqnarray}
in the cycle given in Fig. \ref{fig:cycle}.

The constructed sequence of random walk satisfies the conditions of {\em Lemma} \ref{lem:backwalk}; hence there exists $\widetilde{T}^*$ and $\mu_{\widetilde{T}}\in(0,1)$ such that
\begin{eqnarray}
\mathbb{P}({\rm label}(i_0)={\rm label}(j_0)\,|\,i_{\widetilde{T}},j_{\widetilde{T}}\in\mathscr{C})\geq \mu_{\widetilde{T}} \nonumber
\end{eqnarray}
holds for any $\widetilde{T}\geq \widetilde{T}^*$. Based on this one finally gets that
\begin{eqnarray}
&    & \mathbb{P}(\mathcal{N}(A_{\sigma_{T:k}}, i)\cap \mathcal{N}(A_{\sigma_{T:k}}, j)\neq \emptyset) \nonumber\\
&\geq& \gamma^{(N-|\chi|)q}\mu_{\widetilde{T}} \,\,=\,\, h_T\,\,\in\,\, (0,1),  \nonumber
\end{eqnarray}
which completes the proof. \hfill $\square$

\textbf{Proof of Lemma \ref{lem:scramblingpositive}: } Denote $\{(i_k,j_k)\}_{k=1}^{\frac{N(N-1)}{2}}=\{(i,j): i\neq j,\, i\in V,\, j\in V\}$ as the set of all possible nodes pairs from $V$. Next,
divide $M:1$ into $\frac{(N-1)N}{2}$ subintervals, i.e., $T:1$, $(2T):(T+1)$, $\cdots$, with the length of each being $T$, and assign each pair of indices $I_k=(i_k,j_k)$ to each of the interval $kT:(k-1)T+1$.

For any $0\leq k<\frac{N(N-1)}{2}$, we choose
\begin{eqnarray}
\hat{i}_k \in \mathcal{N}(A_{M:(kT+1)},i_k), \quad
\hat{j}_k \in \mathcal{N}(A_{M:(kT+1)},j_k). \nonumber
\end{eqnarray}
According to {\em Lemma} \ref{lem:sharingprob}, for the pair of indices $\hat{i}_k$ and $\hat{j}_k$, there exist $c_0$ and $\beta$ such that
\begin{eqnarray}
&    & \mathbb{P}(\mathcal{N}(A_{\sigma_{kT:((k-1)T+1)}}, \hat{i}_k)\cap \mathcal{N}(A_{\sigma_{kT:((k-1)T+1)}}, \hat{j}_k)\neq \emptyset) \nonumber\\
&\geq& h_T. \nonumber
\end{eqnarray}
It should be especially noted that if $\mathcal{N}(A_{\sigma_{kT:((k-1)T+1)}}, \hat{i}_k)\cap \mathcal{N}(A_{\sigma_{kT:((k-1)T+1)}}, \hat{j}_k)\neq \emptyset$ holds for any $k\leq \frac{(N-1)N}{2}$, then
\begin{eqnarray}
\mathcal{N}(A_{\sigma_{T:1}}, i)\cap \mathcal{N}(A_{\sigma_{T:1}}, j)\neq \emptyset\nonumber
\end{eqnarray}
for any pair of nodes $i,j\in V$. Based on this, one knows that $A_{\sigma_{M:1}}$ is scrambling with the probability no less than  $h_T^\frac{N(N-1)}{2}$. Note that $\lambda(A_{\sigma_{M:1}})\leq 1-\delta^M$ when $A_{\sigma_{M:1}}$ is scrambling, the proof is hence completed. \hfill $\square$

Based on the above lemmas, we present the  proof of {\em Theorem} \ref{them:1}.

\textbf{Proof of Theorem \ref{them:1}: } Given any sufficiently small $\varepsilon>0$, choose a fixed integer $T\geq T^*$ with $T^*$ given in {\em Lemma} \ref{lem:scramblingpositive}. Based on $T$ choose $q$ such that
$$\varepsilon> (1-\delta^{\frac{N(N-1)}{2}T})^q,$$
where $\delta$ is the minimal positive entry of $A$.

Denote
\begin{eqnarray}
D = \frac{N(N-1)}{2}T, \quad
M = Dq, \quad
d = \left\lfloor\frac{m}{M}\right\rfloor. \nonumber
\end{eqnarray}
Based on $M$ and $d$ one has that
\begin{eqnarray}
&    & \mathbb{P}(\lambda(A_{\sigma_{m:1}})\geq \varepsilon) \nonumber\\
&\leq& \mathbb{P}(\lambda(A_{\sigma_{(dM):1}})\geq \varepsilon) \nonumber\\
&\leq& \mathbb{P}\left(\prod_{r=1}^{d}\lambda(A_{\sigma_{(rM):((r-1)M+1)}})\geq \varepsilon\right)\nonumber\\
&=& 1-\mathbb{P}\left(\prod_{r=1}^{d}\lambda(A_{\sigma((rM):((r-1)M+1))})< \varepsilon\right).\nonumber
\end{eqnarray}

It should be noted that $\prod_{r=1}^{d}\lambda(A_{\sigma((rM):((r-1)M+1))})< \varepsilon$ holds as long as
\begin{eqnarray}
\lambda(A_{\sigma((rM):((r-1)M+1))})< \varepsilon \nonumber
\end{eqnarray}
holds for some $r=1,2,\cdots, d$. In fact, the probability of $\lambda(A_{\sigma((rM):((r-1)M+1))})< \varepsilon$ for some $r\in\{1,2,\cdots, d\}$ is
\begin{eqnarray}
1- \prod_{r=1}^{d}\mathbb{P}(\lambda(A_{\sigma((rM):((r-1)M+1))})\geq \varepsilon). \nonumber
\end{eqnarray}
Hence,
\begin{eqnarray}
\mathbb{P}(\lambda(A_{\sigma_{m:1}})\geq \varepsilon) \leq \prod_{r=1}^{d}\mathbb{P}(\lambda(A_{\sigma((rM):((r-1)M+1))})\geq \varepsilon). \nonumber
\end{eqnarray}

Furthermore, since $M=Dq$, it follows
\begin{eqnarray}
&    & \mathbb{P}\left(\lambda(A_{\sigma((rM):((r-1)M+1))})\geq \varepsilon \right) \nonumber\\
&\leq& 1- \mathbb{P}\left(\prod_{h=1}^q\lambda(A_{\sigma(((r-1)M+hD):((r-1)M+(h-1)D+1))}) \right. \nonumber\\
& &  \left.\qquad\qquad<\varepsilon \frac{}{}\right)\nonumber\\
&\leq& 1- \prod_{h=1}^q\mathbb{P}\left(\lambda(A_{\sigma(((r-1)M+hD):((r-1)M+(h-1)D+1))}) \frac{}{}\right. \nonumber\\
& &  \left.\qquad\qquad<\varepsilon^{\frac{1}{q}} \frac{}{}\right)\nonumber\\
&\leq&
1- \prod_{h=1}^q\mathbb{P}\left(\frac{}{}\lambda(A_{\sigma(((r-1)M+hD):((r-1)M+(h-1)D+1))}) \right. \nonumber\\
&    & \left.\qquad\qquad \leq 1-\delta^D \frac{}{}\right) \nonumber\\
&\leq   & 1- \left(h_T^{\frac{N(N-1)}{2}}\right)^q \nonumber\\
&=   & \gamma^*, \nonumber
\end{eqnarray}
where $\gamma^*\in(0,1)$ and the last inequality follows from {\em Lemma} \ref{lem:scramblingpositive}.

Summarizing the above deduction, it must be true that
\begin{eqnarray}
\mathbb{P}(\lambda(A_{\sigma(m:1)})\geq \varepsilon) < (\gamma^*)^d.\nonumber
\end{eqnarray}
Since $d\rightarrow\infty$ as $m\rightarrow\infty$, one derives that $\lim_{m\rightarrow\infty}$ $\mathbb{P}(\lambda(A_{\sigma(m:1)})\geq \varepsilon)=0$ and the proof is hence completed by directly using {\em Lemma} \ref{lem:scramblingconvergence}. \hfill $\square$

\section{Corollaries of {\em Theorem} \ref{them:1}}

If there exists a node of $\mathbbm{r}(A)$ whose diagonal entry in $A$ is positive, the conditions in {\em Theorem} \ref{them:1} can be greatly simplified.

\begin{coro}\label{coro:0}
The asynchronous DCA (\ref{eq:asynchronous}) generates consensus almost surely if all the following conditions hold:
\begin{itemize}
\item[a)] $\mathcal{G}(A)$ is rooted with one of $\mathbbm{r}(A)$ contains a self-loop, and $\sigma_k\in 2^V$.
\item[b)] There exists $\alpha>0$ such that if $\mathbb{P}(\sigma_{k}\,|\,\sigma_{(k-1):1})\neq 0$, then $\mathbb{P}(\sigma_{k}\,|\,\sigma_{(k-1):1})\geq \alpha$.
\item[c)] There exists $q>0$ such that
\begin{eqnarray}
\bigcup_{w=0}^{q-1}\left(\bigcup_{\sigma\in\mathscr{I}_{\sigma_{(k-1):1}(w)}}\sigma\right)= V, \quad\forall\, k\geq 1, \nonumber
\end{eqnarray}
holds for any given values of $\sigma_{(k-1):1}$, where
\begin{eqnarray}
\mathscr{I}_{\sigma_{(k-1):1}}(w)=\{\sigma: \mathbb{P}(\sigma_{k+w}=\sigma\,|\,\sigma_{(k-1):1})\neq 0\}. \nonumber
\end{eqnarray}
\end{itemize}
\end{coro}

Conditions b) and c) in {\em Corollary} \ref{coro:0} cover Markovian switching sequence as a special case, which is similar to the case in \cite{SIAM-Matei}. Furthermore, condition b) in {\em Corollary} \ref{coro:0} (also, in {\em Theorem} \ref{them:1}) is difficult to be relaxed, which can be seen from the following example.

\begin{exam}
Let $\sigma_k\in V$ and consider the stochastic matrix
\begin{eqnarray}
A = \left(\begin{array}{ccc}
0.5 & 0.5 & 0\\
0 & 0 & 1\\
1 & 0 & 0
\end{array}\right) \nonumber
\end{eqnarray}
and the following Markovian transition probability from $\sigma_k$ to $\sigma_{k+1}$:
\begin{eqnarray}
M_k = \left(\begin{array}{ccc}
1-\frac{1}{k} & 0 & 1\\
\frac{1}{k} & 0 & 0\\
0 & 1 & 0
\end{array}\right). \nonumber
\end{eqnarray}
If $\sigma_1=1$, the corresponding product of stochastic matrices is
\begin{eqnarray}
& & \cdots A_{\sigma_{k}}\cdots A_{\sigma_{2}}A_{\sigma_{1}} \nonumber\\
&=& \cdots A_{3}A_{2}A_{1}^{k_p}\cdots A_{3}A_{2}A_{1}^{k_2}A_{3}A_{2}A_{1}^{k_1}, \nonumber
\end{eqnarray}
where $\lim_{p\rightarrow\infty}k_p\overset{\text{a.s.}}{=}+\infty$. One can verify that for any $k\geq 1$, it holds
\begin{eqnarray}
A_1^k &=&
\left(\begin{array}{ccc}
\frac{1}{2^k} & 1-\frac{1}{2^k} & 0\\
0 & 1 & 0\\
0 & 0 & 1
\end{array}\right), \nonumber\\
A_3A_2A_1^k &=&
\left(\begin{array}{ccc}
\frac{1}{2^k} & 1-\frac{1}{2^k} & 0\\
0 & 0 & 1\\
\frac{1}{2^k} & 1-\frac{1}{2^k} & 0
\end{array}\right). \nonumber
\end{eqnarray}
We further define
\begin{eqnarray}
w_1 = (0, 1, 0)^T, \quad
w_2 = (1, 0, 1)^T, \nonumber
\end{eqnarray}
and obtain that
\begin{eqnarray}
A_3A_2A_1^k w_1 &=& (1-\frac{1}{2^k})w_2, \nonumber\\
A_3A_2A_1^k w_2 &=& (1-\frac{1}{2^k})w_1 + \frac{1}{2^k}\mathbbold{1}, \nonumber
\end{eqnarray}
based on which one derives that
\begin{eqnarray}
& & A_{3}A_{2}A_{1}^{k_p}\cdots A_{3}A_{2}A_{1}^{k_2}A_{3}A_{2}A_{1}^{k_1} \nonumber\\
&=& \left\{\begin{array}{ll}
\prod_{\gamma=1}^p (1-\frac{1}{2^{k_\gamma}}) w_1 + \sum_{\gamma=1}^{p/2}\frac{1}{2^{k_\gamma}}\mathbbold{1}, &\!\!\!\! \mbox{ if } p \mbox{ is even}, \\
\prod_{\gamma=1}^p (1-\frac{1}{2^{k_\gamma}}) w_2 + \sum_{\gamma=1}^{(p-1)/2}\frac{1}{2^{k_\gamma}}\mathbbold{1}, &\!\!\!\! \mbox{ if } p \mbox{ is odd}.
\end{array}\right.\nonumber
\end{eqnarray}
Since $w_1\neq w_2$ and $\prod_{\gamma=1}^p (1-\frac{1}{2^{k_\gamma}})$ converges to a nonzero value almost surely as $p$ tends to infinity, one knows that $A_{\sigma_{k:1}}$ does not converge to consensus almost surely as $k$ goes to infinity.
\end{exam}

If the random variables $\{\Xi_k\}_{k\geq 1}$ are independent, then condition c) in {\em Theorem} \ref{them:1} can be removed and one obtains the following corollary.

\begin{coro}\label{coro:1}
The asynchronous DCA (\ref{eq:asynchronous}) generates consensus almost surely if all the following conditions hold:
\begin{itemize}
\item[a)] $\mathcal{G}(A)$ is rooted and $\sigma_k\in 2^V$.
\item[b)] $\{\Xi_k\}_{k\geq 1}$ are independent and there exists $\alpha>0$ such that if $\mathbb{P}(\sigma_{k})\neq 0$, then $\mathbb{P}(\sigma_{k})\geq \alpha$.
\item[c)] There exists $q>0$ such that
\begin{eqnarray}
\bigcup_{\tau=k}^{k+q-1}\left(\bigcup_{\sigma\in\mathscr{I}_\tau}\sigma\right)= V, \quad\forall\, k\geq 1, \nonumber
\end{eqnarray}
where
$$\mathscr{I}_{k}= \{\sigma: \mathbb{P}(\sigma_{k}=\sigma)\neq 0\}.$$
\item[d)] There exists a strongly connected component $\chi\subseteq\mathbbm{r}(A)$ such that for any $j\in \chi$, there is $\mathscr{I}_{k}^{j}\neq\emptyset$ and
    \begin{eqnarray}
    \chi\bigcap\left(\bigcap_{\sigma\in \mathscr{I}_{k}^{j} } \sigma\right) = j, \quad \forall\, k\geq 1,  \nonumber
    \end{eqnarray}
    where
    $$\mathscr{I}_{k}^{j}= \{\sigma: \sigma \in \mathscr{I}_k \mbox{ and } j\in \sigma\}.$$
\end{itemize}
\end{coro}

Condition d) in {\em Corollary} \ref{coro:1} (also, in {\em Theorem} \ref{them:1}) is difficult to be relaxed. In the following example, we find that consensus behavior does not happen if we remove condition d) in {\em Corollary} \ref{coro:1}.

\begin{exam}
Let
\begin{eqnarray}
A =\left(\begin{array}{cccc}
0 &0 &0 &1\\
1 &0 &0 &0\\
0 &1 &0 &0\\
0 &0 &1 &0
\end{array}\right) \nonumber
\end{eqnarray}
and choose
\begin{eqnarray}
\mathbb{P}(\sigma_{k}) =
\left\{\begin{array}{ll}
1,  &  \mbox{ if } k\equiv 1\,({\rm mod}\,\, 4) \, \mbox{ and } \sigma_k=\{1,3\},\\
0.5,  &  \mbox{ if } k\equiv 2\,({\rm mod}\,\, 4) \, \mbox{ and } \sigma_k=\{1\},\{3\},\\
1,  &  \mbox{ if } k\equiv 3\,({\rm mod}\,\, 4) \, \mbox{ and } \sigma_k=\{2,4\},\\
0.5,  &  \mbox{ if } k\equiv 0\,({\rm mod}\,\, 4) \, \mbox{ and } \sigma_k=\{2\},\{4\}.
\end{array}\right. \nonumber
\end{eqnarray}
Then one can verify that condition c) holds for $q=4$, and
\begin{eqnarray}
\mathscr{I}_k^{1}=
\{1,3\},  \quad \mbox{if}\,\,\, k\equiv 1({\rm mod}\,\, 4). \nonumber
\end{eqnarray}
However, since $\mathbbm{r}(A)=\{1,2,3,4\}$ and any proper subset of $\mathbbm{r}(A)$ is not a strongly connected component, one sets $\chi=\mathbbm{r}(A)$ and knows that
\begin{eqnarray}
    \chi\bigcap\left(\bigcap_{\sigma\in \mathscr{I}_{k}^{1} } \sigma\right) = \{1,3\}\neq \{1\},\quad \mbox{if}\,\,\, k\equiv 1({\rm mod}\,\, 4), \nonumber
\end{eqnarray}
which indicates condition d) is violated. Moreover,
\begin{eqnarray}
A_{\sigma_{k:1}}= \left\{
\begin{array}{ll}
A_{\{1,3\}}(A_{\{2,4\}}A_{\{1,3\}})^\frac{k-1}{4},  & \mbox{ if } k\equiv 1\,({\rm mod}\,\, 4),\\
A_{\{1,3\}}(A_{\{2,4\}}A_{\{1,3\}})^\frac{k-2}{4},  & \mbox{ if } k\equiv 2\,({\rm mod}\,\, 4),\\
(A_{\{2,4\}}A_{\{1,3\}})^\frac{k+1}{4},  & \mbox{ if } k\equiv 3\,({\rm mod}\,\, 4),\\
(A_{\{2,4\}}A_{\{1,3\}})^\frac{k}{4},  & \mbox{ if } k\equiv 0\,({\rm mod}\,\, 4),
\end{array}
\right.\nonumber
\end{eqnarray}
since
\begin{eqnarray}
A_{\{2,4\}}A_{\{1,3\}}=
\left(\begin{array}{cccc}
0& 0& 0& 1\\
0& 0& 0& 1\\
0& 1& 0& 0\\
0& 1& 0& 0
\end{array}\right) \nonumber
\end{eqnarray}
is not SIA, one knows $A_{\sigma_{k:1}}$ does not converge.
\end{exam}

The conditions in {\em Corollary} \ref{coro:1} does not imply the {\em\color{red} strongly aperiodic condition}. In \cite{I3EAC-Touri-Nedic}, Touri and Nedi\'{c} proposed a strongly aperiodic condition for convergence of products of random stochastic matrices, which means there exists $\hat{\gamma}\in(0,1)$ such that
\begin{eqnarray}
\mathbb{E}(A_{\sigma_k}(i,i)A_{\sigma_k}(i,j))\geq \hat{\gamma} \mathbb{E}(A_{\sigma_k}(i,j)) \label{eq:strongergodic}
\end{eqnarray}
for any $i,j\in V$ ($i\neq j$). In the following example, we will show that (\ref{eq:strongergodic}) may not hold even if the conditions of {\em Corollary} \ref{coro:1} are all satisfied.

\begin{exam}
Suppose $\{\Xi_k\}_{k\geq 1}$ are i.i.d, $V=\{1,2,3,4\}$, and $A$ is given in (\ref{eq:eg4}). When we set $\mathscr{I}_k=\{1,2,3,4\}$ and $\mathbb{P}(\sigma_k=j)=\frac{1}{4}$ for any $j\in \mathscr{I}_k$, the correctness of the $4$ conditions in {\em Corollary} \ref{coro:1} can be easily verified. Using simple calculation, one has that
\begin{eqnarray}
\mathbb{E}(A_{\sigma_k}(1,1)A_{\sigma_k}(1,2)) = \frac{1}{4}\cdot 0 + \frac{1}{4}\cdot 0 + \frac{1}{4}\cdot 0 + \frac{1}{4}\cdot 0 = 0. \nonumber
\end{eqnarray}
However,
\begin{eqnarray}
\mathbb{E}(A_{\sigma_k}(1,2)) = \frac{1}{4}\cdot 1 + \frac{1}{4}\cdot 0 + \frac{1}{4}\cdot 0 + \frac{1}{4}\cdot 0 = \frac{1}{4} ,\nonumber
\end{eqnarray}
which contradicts equation (\ref{eq:strongergodic}).
\end{exam}

When each $\sigma_k$ ($k\geq 1$) takes values in $V$, {\em Theorem} \ref{them:1} can be simplified in the following form.

\begin{coro}\label{coro:2}
The asynchronous DCA (\ref{eq:asynchronous}) generates consensus almost surely if all the following conditions hold:
\begin{itemize}
\item[a)] $\mathcal{G}(A)$ is rooted and $\sigma_k\in V$.
\item[b)] There exists $\alpha>0$ such that if $\mathbb{P}(\sigma_{k}\,|\,\sigma_{(k-1):1})\neq 0$, then $\mathbb{P}(\sigma_{k}\,|\,\sigma_{(k-1):1})\geq \alpha$.
\item[c)] For any given historic values of $\sigma_{(k-1):1}$, the set
    $$\mathscr{I}_{\sigma_{(k-1):1}}= \{\sigma: \mathbb{P}(\sigma_{k}=\sigma\,|\,\sigma_{(k-1):1})\neq 0\}$$
    only depends on $k$ but not the historic values $\sigma_{(k-1):1}$, i.e, there exists $\mathscr{I}_{k}$ such that
    $$\mathscr{I}_{\sigma_{(k-1):1}}=\mathscr{I}_{k},\quad \forall\, \sigma_{(k-1):1}\in \underbrace{V\times V\cdots \times V}_{k-1\,\, \mbox{\footnotesize times}}. $$
\item[d)] There exists $q>0$ such that
\begin{eqnarray}
\bigcup_{\tau=k}^{k+q-1}\mathscr{I}_{\tau}= V, \quad\forall\, k\geq 1, \nonumber
\end{eqnarray}
where
$$\mathscr{I}_{k}= \{\sigma: \mathbb{P}(\sigma_{k}=\sigma\,|\,\sigma_{k:1})\neq 0\}.$$
\item[e)] There exists a strongly connected component $\chi\subseteq \mathbbm{r}(A)$ such that: $\chi\subseteq \mathscr{I}_{k}$ holds for any $k\geq 1$.
\end{itemize}
\end{coro}

Condition c) in {\em Corollary} \ref{coro:2} (also, in {\em Theorem} \ref{them:1}) plays a very critical role for consensus and is difficult to be relaxed, which can be observed from the following example.

\begin{exam}
Consider the following stochastic matrix
\begin{eqnarray}
A = \left(\begin{array}{ccc}
0 & 0 & 1\\
1 & 0 & 0\\
0 & 1 & 0
\end{array}\right) \nonumber
\end{eqnarray}
and use the following ergodic Markovian chain $M$ (also has positive diagonal entries) to generate the switching signal $\sigma_k$ in DCA (\ref{eq:asynchronous}):
\begin{eqnarray}
M = \left(\begin{array}{ccc}
0.5 & 0.5 & 0\\
0 & 0.5 & 0.5\\
0.5 & 0 & 0.5
\end{array}\right). \nonumber
\end{eqnarray}
One can verify that $\mathscr{I}_{\sigma_1=1}=\{1,3\}$, $\mathscr{I}_{\sigma_1=2}=\{1,2\}$, and $\mathscr{I}_{\sigma_1=3}=\{2,3\}$, which violates condition c) of {\em Corollary} \ref{coro:2}. When $\sigma_1=3$, the generated products $\prod_{k=1}^{\infty}A_{\sigma_k}$ has the form of
\begin{eqnarray}
& & \cdots A_{\sigma_k}A_{\sigma_{k-1}}\cdots A_{\sigma_1} \nonumber\\
&=& \cdots A_1^{k_{3p}}A_2^{k_{3p-1}}A_3^{k_{3p-2}}\cdots A_1^{k_6}A_2^{k_5}A_3^{k_4}A_1^{k_3}A_2^{k_2}A_3^{k_1}, \nonumber
\end{eqnarray}
where $k_\tau\geq 1$ for any $\tau\geq 1$. By using simple calculation, one knows that
\begin{eqnarray}
A_1^k &=& A_1, \quad k\geq 1, \nonumber\\
A_2^k &=& A_2, \quad k\geq 1, \nonumber\\
A_3^k &=& A_3, \quad k\geq 1, \nonumber\\
A_1A_2A_3 &=& \bar{A} \,\,\,= \,\,\,\left(\begin{array}{ccc}
0 & 1 & 0\\
1 & 0 & 0\\
0 & 1 & 0
\end{array}\right), \nonumber
\end{eqnarray}
where $\bar{A}$ is not SIA. Hence,
\begin{eqnarray}
& &\cdots A_1^{k_{3p}}A_2^{k_{3p-1}}A_3^{k_{3p-2}}\cdots A_1^{k_6}A_2^{k_5}A_3^{k_4}A_1^{k_3}A_2^{k_2}A_3^{k_1} \nonumber\\
&=& \cdots(A_1A_2A_3)(A _1A_2A_3)(A_1A_2A_3) \,\,\,= \,\,\, \cdots\bar{A}\bar{A}\bar{A} \nonumber
\end{eqnarray}
does not converge to consensus.
\end{exam}

If we require $j\in \mathscr{I}_k$ for any $k\geq 1$ and $j\in V$ in {\em Theorem} \ref{them:1}, one further derives the following result.

\begin{coro}\label{coro:3}
The asynchronous DCA (\ref{eq:asynchronous}) generates consensus almost surely if all the following conditions hold:
\begin{itemize}
\item[a)] $\mathcal{G}(A)$ is rooted and $\sigma_k\in 2^V$.
\item[b)] There exists $\alpha\in(0,1)$ such that if $\mathbb{P}(\sigma_{k}\,|\,\sigma_{(k-1):1})\neq 0$, then $\mathbb{P}(\sigma_{k}\,|\,\sigma_{(k-1):1})\geq \alpha$.
\item[c)] $\mathbb{P}(\sigma_{k}=j\,|\,\sigma_{(k-1):1})\neq 0, \quad \forall j\in V,\, k\geq 1$.
\end{itemize}
\end{coro}

The proof of {\em Corollary} \ref{coro:1}-\ref{coro:3} can be derived similarly to that of {\em Theorem} \ref{them:1} and hence omitted here.

\end{document}